\title{
Simple Proofs of Classical Theorems in Discrete Geometry via
the Guth--Katz Polynomial Partitioning Technique\thanks{%
        Work by Haim Kaplan has been supported
        by Grant 2006/204 from the U.S.-Israel Binational Science Foundation,
        and by grant 822/10 from the Israel Science Fund.
        Work by Micha Sharir has been supported
        by NSF Grant CCF-08-30272,
        by Grant 2006/194 from the U.S.-Israel Binational Science Foundation,
        by grant 338/09 from the Israel Science Fund,
        and by the Hermann Minkowski--MINERVA Center for Geometry at Tel Aviv
        University.}}

\documentclass[11pt]{article}

\usepackage{latexsym,color,graphicx,amsmath,amssymb}
\usepackage{epsfig}

\usepackage{fullpage}

\newtheorem{theorem}{Theorem}[section]
\newtheorem{lemma}[theorem]{Lemma}

\newcommand{\Deg}{D}  

\author{
{\sc Haim Kaplan}
\\
    {\footnotesize School of Computer Science, }\\[-1.5mm]
    {\footnotesize Tel Aviv University, }\\[-1.5mm]
    {\footnotesize Tel~Aviv 69978, Israel}
\and {\sc Ji\v{r}\'{\i} Matou\v{s}ek}
\\
   {\footnotesize Department of Applied Mathematics and}\\[-1.5mm]
   {\footnotesize Institute of Theoretical Computer Science (ITI)}\\[-1.5mm]
   {\footnotesize  Charles University, Malostransk\'{e} n\'{a}m. 25}\\[-1.5mm]
   {\footnotesize  118~00~~Praha~1, Czech Republic, and}\\
   {\footnotesize    Institute of  Theoretical Computer Science}\\[-1.5mm]
   {\footnotesize    ETH Zurich, 8092 Zurich, Switzerland}
\and {\sc Micha Sharir}
\\
    {\footnotesize School of Computer Science, }\\[-1.5mm]
    {\footnotesize Tel Aviv University, }\\[-1.5mm]
    {\footnotesize Tel~Aviv 69978, Israel, and}\\
    {\footnotesize Courant Institute of Mathematical Sciences, }\\[-1.5mm]
    {\footnotesize New York University, }\\[-1.5mm]
    {\footnotesize New York, NY~~10012,~USA}
}

\newcommand{\heading}[1]{\vspace{1ex}\par\noindent{\bf #1}}
\newcommand\notion[1]{\emph{#1}}

\makeatletter
\newcommand{\ProofEndBox}{{\ifhmode\unskip\nobreak\hfil\penalty50 \else
          \leavevmode\fi\quad\vadjust{}\nobreak\hfill$\Box$
            \finalhyphendemerits=0 \par}}
\makeatother
\newcommand{\proofend}{\ProofEndBox\smallskip}

\newcommand{\R}{{\mathbb{R}}}

\newcommand\PP{\mathcal{P}}

\def\:{\colon}

\long\def\onefigure#1#2{
\begin{figure*}[tbp]
\begin{center}
#1
\end{center}
\caption{#2}
\end{figure*}
}

\newcommand{\labepsfig}[2]  
{\onefigure{\mbox{\epsfig{file=pst-#1.eps}}}{\label{f:#1}
{\small\sf #2}} }

\newcommand{\labepsfigw}[3]  
{\onefigure{\mbox{\epsfig{file=#1.eps,width=#2}}}{\label{f:#1}
{\small\sf #3}} }

\begin{document}

\maketitle

\begin{abstract}
Recently Guth and Katz \cite{GK2} invented, as a step in their
nearly complete solution of Erd\H{o}s's distinct distances problem,
a new method for partitioning finite point sets in $\R^d$, based on
the Stone--Tukey polynomial ham-sandwich theorem. We apply this
method to obtain new and simple proofs of two well known results:
the Szemer\'edi--Trotter theorem on incidences of points and lines,
and the existence of spanning trees with low crossing numbers. Since
we consider these proofs particularly suitable for teaching, we aim
at self-contained, expository treatment. We also mention some
generalizations and extensions, such as the Pach--Sharir bound on
the number of incidences with algebraic curves of bounded degree.
\end{abstract}

\section{Introduction}

A dramatic breakthrough in discrete geometry took place in
November 2010, when Guth and Katz \cite{GK2} completed a
project of Elekes, exposed in \cite{ES}, and established a
nearly complete solution of Erd\H{o}s's  distinct distances problem
\cite{E46}, originally posed in 1946.

In one of the main steps of their analysis, they apply the
\emph{polynomial ham-sandwich theorem} of Stone and Tukey~\cite{ST}
to obtain a partition of a finite point set $P$ in $\R^d$ with
certain favorable properties, detailed in Section~\ref{s:discuss-partit} below.
The partition is
effected by what we call an \emph{$r$-partitioning polynomial}.
The removal of the zero set $Z$ of the polynomial
partitions space into connected components, each
 containing at most $|P|/r$ points of~$P$. A key feature of the
construction is that the degree of the polynomial achieving this
need not be too high, only $O(r^{1/d})$, and thus the interaction
of other objects, such as lines or hyperplanes, with $Z$ is under control in
some sense.

In this paper we apply  partitioning polynomials in
several classical problems of discrete geometry, mostly
planar ones, and we provide new and simple proofs of
some well known results.


\paragraph{Incidences.}
For a finite set $P\subset\R^2$ and a finite set $L$ of lines
in $\R^2$, let $I(P,L)$ denote the number of
\notion{incidences} of $P$ and $L$, i.e., of pairs
$(p,\ell)$ with $p\in P$, $\ell\in L$, and $p\in\ell$.

The following fundamental result was first proved by Szemer\'edi and Trotter
in 1983, in response to a problem of Erd\H{o}s \cite{E46}.

\begin{theorem}[Szemer\'edi and Trotter \cite{SzT}] \label{thm:SZ}
$I(P,L)=O(m^{2/3}n^{2/3}+m+n)$ for every set $P$ of $m$
distinct points in the plane and every set $L$ of $n$ distinct
lines.
\end{theorem}

We remark that the bound in the theorem is tight in the worst case
for all $m,n$ (see \cite{E46}, \cite{E01} for original
sources or \cite{Ma02} for a presentation).

A simpler proof of the Szemer\'edi--Trotter theorem,
based on cuttings, was given by Clarkson et al.~\cite{CEGSW} in
1990, and in 1997 Sz\'ekely~\cite{s-cnhep-97} found a beautiful
and elegant proof, based on the crossing lemma for graphs embedded
in the plane (also see, e.g., \cite{Ma02}).

In Section~\ref{s:SzTproof}
we present an alternative proof based on polynomial partitions,
hoping that the reader will find it equally simple. We also
believe that the new proof is  suitable for teaching purposes, so
our goal is to make the exposition as elementary and self-contained
as possible. For this we also give proofs of several well known and
 basic facts about multivariate polynomials. The only major ingredient of
the analysis which we do not prove is the classical ham-sandwich
theorem, which we use as a black box (see, e.g.,
\cite{Ma03} for an exposition).

The Szemer\'edi--Trotter theorem has led to an extensive study of
incidences of points and curves in the plane and of points and
surfaces in higher dimensions. A survey of the topic can be found in
Pach and Sharir~\cite{incsurv}. In particular, the following
theorem on incidences between points and planar curves has been established:

\begin{theorem}[Pach and Sharir~\cite{ps-nibpc-98}] \label{szt:algeb}
Let $P$ be a set of $m$ points and let $\Gamma$ be a set of $n$
simple curves, all lying in the plane. If no more than $C_1$ curves
of $\Gamma$ pass through any $k$ given points, and every pair of
curves of $\Gamma$ intersect in at most $C_2$ points, then
$$
I(P,\Gamma) = O\left(m^{k/(2k-1)}n^{(2k-2)/(2k-1)}+m+n \right),
$$
with an appropriate constant of proportionality that depends on $k,C_1,C2$.
\end{theorem}

A weaker version of this result, where the the curves
in $\Gamma$ are assumed to be algebraic and to belong to a family
parameterized by $k$ real parameters, was obtained earlier,
also by Pach and Sharir \cite{ps-raprp-92}
(special cases of this result, e.g., for incidences of points
and circles, were obtained even earlier by Clarkson et al.~\cite{CEGSW}).

In Section~\ref{s:curves}, we give a simple proof of
a version of Theorem~\ref{szt:algeb},
with the additional
assumption that $\Gamma$ consists of algebraic curves of
degree bounded by a constant.

\paragraph{Spanning trees with low crossing number.}
Let $P$ be a finite set of points in $\R^2$. A \emph{(geometric)
graph} on $P$ is a graph $G$ with vertex set $P$
whose edges are realized
as straight segments connecting the respective end-vertices.
The \emph{crossing number} of $G$ is
the maximum number of edges that can be intersected simultaneously
by a line not passing through any point of~$P$.\footnote{The condition
of avoiding the points of $P$ is important; for example,
if all of the points of $P$ are collinear, then the line
containing $P$ necessarily intersects all edges.}
We will consider \emph{geometric spanning trees} on $P$,
i.e., acyclic connected geometric graphs on~$P$.

The following result has been established in the late
1980s by Welzl \cite{w-pttco-88} and by Chazelle and Welzl \cite{CW};
also see  Welzl \cite{We}.

\begin{theorem}[Welzl \cite{w-pttco-88}, Chazelle and Welzl \cite{CW}]
\label{t:sptree}
Every set of $n$ points in the plane has a geometric spanning tree
with crossing number $O(\sqrt n\,)$.
\end{theorem}

The bound in the theorem is tight up to a multiplicative constant,
as the example of a $\sqrt n\times\sqrt n$ grid shows.
Spanning trees with low crossing number have many applications in
discrete and computational geometry, including
range searching \cite{CW}, the design of other geometric algorithms
(see, e.g., \cite{Asano-al-barriers}), discrepancy theory \cite{MWW},
and approximation \cite{m-zono}.

The original proof of Theorem~\ref{t:sptree} constructs
the tree iteratively, through a process called
\emph{iterative reweighting}.  In each step several new
edges are added, and these are selected
using a packing argument with balls
in a line arrangement (or, alternatively, using
a so-called cutting).
An alternative proof, replacing
iterative reweighting with linear programming duality,
was recently given by Har-Peled \cite{Sariel}.

In Section~\ref{s:pl-sp}
we present a new and simple proof of Theorem~\ref{t:sptree}
via polynomial partitions.

Chazelle and Welzl \cite{CW} established their result on spanning trees
with low crossing number in a very general setting,
where the points do not lie in the plane, but rather in the
ground set of an arbitrary set system $\mathcal F$.
The bound on the crossing number is then expressed in terms
of the \emph{dual shatter function} of~$\mathcal F$.

At present it seems that the approach with polynomial partitions is
not suitable for this level of generality. However, some
generalizations are possible. First, we have verified that
Theorem~\ref{t:sptree} can be extended to the case where the
crossing number is taken with respect to a family of algebraic
curves of degree bounded by a constant, but
we will not pursue this in this paper.

Second, one can also prove a $d$-dimensional
generalization of Theorem~\ref{t:sptree}, and this we do
in Section~\ref{s:high-sp}.
Here we are given a set $P$ of $n$ points in $\R^d$, and
consider spanning trees of $P$, which we embed into $\R^d$ by
drawing their edges as straight segments connecting the respective
end-vertices, as in the plane. The crossing number of such a tree is
the maximum number of its edges that are crossed by a hyperplane not
passing through any point of $P$.
According to \cite{w-pttco-88}, \cite{CW}, for every
$n$-point set in $\R^d$ there exists a
straight-edge spanning tree with crossing number
$O\left(n^{1-1/d}\right)$.

We re-prove this fact
using polynomial partitions, similar to the planar case.
However, the proof is
more involved in higher dimensions. Informally, the partition
distributes the input points evenly among the resulting cells,
except that some (in the worst case even all) of the points may
lie on the zero set $Z$ of the partitioning polynomial, and
therefore not belong to any of the subsets.

We avoid this situation using a perturbation argument.
This works for the spanning tree construction because there we may assume
general position of the input points.
For incidence problems this
assumption cannot be made, and other techniques are needed to handle
the points on~$Z$. We intend to investigate  alternative approaches
to handling points on $Z$ in a subsequent paper.

\section{Review of tools}

\subsection{Preliminaries on polynomials}

Here we recall some  standard facts about polynomials. The proofs
are given for didactic purposes, and can be skipped by more
experienced readers.

Since most of the problems that we study here are planar, we will
consider mostly bivariate polynomials
$f=f(x,y)=\sum_{i,j} a_{ij} x^i y^j\in\R[x,y]$, but the analysis can
easily be extended to $d$-variate polynomials in $\R^d$.
The \notion{degree} of $f$ is $\deg(f)=\max\{i+j \mid a_{ij}\ne 0\}$.
Let $Z(f)=\{(x,y)\in\R^2 \mid f(x,y)=0\}$ denote
the \notion{zero set} of $f$.

\begin{lemma}\label{l:atmost}
If $\ell$ is a line in $\R^2$ and $f\in \R[x,y]$
is of degree at most $\Deg$, then either $\ell\subseteq Z(f)$,
or $|\ell\cap Z(f)|\le \Deg$.
\end{lemma}

\heading{Proof.}
Writing $\ell$ in parametric form
$\{(u_1t+v_1,u_2t+v_2) \mid t\in\R\}$, we get that the points of $\ell
\cap Z(f)$ are roots of the univariate polynomial
$g(t):=f(u_1t+v_1,u_2t+v_2)$, which is of degree at most~$\Deg$.
Thus, either $g$ is identically $0$, or it has at most $\Deg$ roots.
\proofend

\begin{lemma}\label{l:contfewlines}
If $f\in \R[x,y]$ is nonzero and of degree at most $\Deg$,
then $Z(f)$ contains at most $\Deg$ distinct lines.
\end{lemma}

\heading{Proof.}
We need to know  that a nonzero bivariate polynomial
(i.e., with at least one nonzero coefficient) does not vanish
on all of $\R^2$. (Readers who do not consider this
a sufficiently standard fact are welcome to work out a quick proof.)

Now we fix a point $p\in\R^2$ not belonging to $Z(f)$.
Let us suppose that $Z(f)$ contains lines $\ell_1,\ldots, \ell_k$.
We choose another line $\ell$ passing through $p$ that is
not parallel to any $\ell_i$ and not passing
through any of the intersections $\ell_i\cap\ell_j$.
(Such an $\ell$ exists since only finitely many directions
need to be avoided.) Then
$\ell$ is not contained in $Z(f)$ and it
has $k$ intersections with $\bigcup_{i=1}^k\ell_i$.
Lemma~\ref{l:atmost} yields $k\le \Deg$.
\proofend

In the proof of Theorem~\ref{t:sptree} (spanning
trees with low crossing number), we will also need
the following result.

\begin{theorem}[Harnack's curve theorem \cite{Harnack}]\label{t:harnack}
Let $f\in\R[x,y]$ be a bivariate polynomial of degree $\Deg$.
Then the number of (arcwise) connected components of $Z(f)$
is at most $1+{\Deg-1\choose 2}$. 
The bound is tight in the worst case.
\end{theorem}

For our application, we actually do not need the precise bound
in Harnack's theorem; it suffices to know that the number
of components is at most~$O(\Deg^2)$. For the sake of completeness,
we provide a short proof of an almost tight bound.

First we recall, without proof, another basic result in algebraic geometry;
see, e.g., \cite{BPR,CLO1,CLO2}.

\begin{theorem}[B\'ezout's theorem]\label{t:bezout}
Let $f,g\in\R[x,y]$ be two bivariate polynomials of degrees
$\Deg_f$ and $\Deg_g$, respectively. (a) If the system $f=g=0$ has
finitely many solutions, then their number is at most $\Deg_f
\Deg_g$. (b) If the system $f=g=0$ has infinitely many solutions,
then $f$ and $g$ have a nontrivial common factor.
\end{theorem}

For a proof of Theorem~\ref{t:harnack}, we choose a generic direction,
and assume, without loss of generality, that it is the
$x$-direction. We may assume that $f$ is square-free, because
eliminating repeated factors of $f$ does not change its zero set.

Every bounded component of $Z(f)$ has at least two extreme
points in the $x$-direction (that is, its leftmost and rightmost
points).
 Such an extreme point has to satisfy $f=f_y=0$,
where $f_y$ is the partial derivative of $f$ with respect to~$y$.

Since $f$ is square-free, $f$ and $f_y$  have no common
factor,\footnote{Assume by induction that this is true for
polynomials of degree smaller than $\Deg$, and let $f$ be a square-free
polynomial of degree $\Deg$. Assume that $f=h\cdot g$ and $f_y = h
\cdot k$ for some polynomials $h$, $g$, and  $k$, where $h$ is
not a constant. Then $f_y = h_y\cdot g + g_y \cdot h = h \cdot
k$. So $h$ divides $h_y\cdot g$. By induction, $h$ and $h_y$ have no
common factors, and so $h$ divides $g$, contradicting our assumption that
$f$ is square-free.} and so by Theorem~\ref{t:bezout} the system
$f=f_y = 0$ has at most $\Deg(\Deg-1)$ solutions. Every bounded
component consumes at least two of these critical points, and hence  the
number of bounded components is at most $\frac12 \Deg(\Deg-1)$.

If $B$ is a sufficiently large number, then (again, assuming generic
directions of the coordinate axes) every unbounded
component of $Z(f)$ meets (at least) one of the two lines
$x=+B$ and $x=-B$. Thus, there are at most $2\Deg$ unbounded
components, and in total we get a bound of $\frac12 \Deg(\Deg+1)$
on all components.


\subsection{The polynomial ham-sandwich theorem}\label{s:veronese}

Here we review the \notion{polynomial ham-sandwich theorem}
of Stone and Tukey \cite{ST}, the key tool used by
Guth and Katz in constructing their partitioning polynomials.

We assume the standard
\notion{ham-sandwich theorem} in the following discrete
version: \emph{Every $d$ finite sets $A_1,\ldots,A_d\subset \R^d$
can be simultaneously bisected by a hyperplane}. Here a hyperplane $h$
\notion{bisects} a finite set $A$ if neither of the two open  halfspaces
bounded by $h$ contains more than $\lfloor|A|/2\rfloor$ points of $A$.

From this, it is easy to derive the polynomial ham-sandwich theorem,
which we state for bivariate polynomials.

\begin{theorem} \label{t:phsc}
Let $A_1,\ldots,A_s\subseteq \R^2$ be finite sets, and let $\Deg$ be
an integer such that ${\Deg+2\choose 2}-1\ge s$. Then there exists a
nonzero polynomial $f\in \R[x,y]$ of degree at most $\Deg$ that
simultaneously bisects all the sets $A_i$, where ``$f$ bisects
$A_i$'' means that $f>0$ in at most $\lfloor|A_i|/2\rfloor$ points
of $A_i$ and $f<0$ in at most $\lfloor|A_i|/2\rfloor$ points of
$A_i$.
\end{theorem}

\heading{Proof.}
We note that ${\Deg+2\choose 2}$ is the number of monomials in
a bivariate polynomial of degree $\Deg$,
or in other words, the number of pairs $(i,j)$ of nonnegative
integers with $i+j\le \Deg$. We set $k:={\Deg+2\choose 2}-1$,
and we let $\Phi\:\R^2\to \R^k$ denote the \notion{Veronese map},
given by
$$
\Phi(x,y) := \left( x^iy^j\right)_{(i,j) \mid 1\le i+j\le \Deg}\in \R^k.
$$
(We think of the coordinates in $\R^k$ as indexed by
pairs $(i,j)$ with $1\le i+j\le \Deg$.)

Assuming, as we may, that $s=k$,  we set $A'_{i}:=\Phi(A_i)$,
$i=1,2,\ldots,k$, and we let $h$ be a hyperplane simultaneously
bisecting $A'_1,\ldots,A'_k$. Then $h$ has an equation of the form
$a_{00}+\sum_{i,j} a_{ij} z_{ij}=0$, where $(z_{ij})_{(i,j) \mid
1\le i+j\le d}$ are the coordinates in $\R^k$. It is easy to check
that $f(x,y):=\sum_{i,j} a_{ij}x^iy^j$ is the desired polynomial
(where here the sum includes $a_{00}$ too).
\proofend

\subsection{Partitioning polynomials}\label{s:discuss-partit}

In this section we recall the construction of Guth and Katz \cite{GK2},
specialized to the planar setting
(our formulation is
 slightly different from theirs). We also (informally) compare it to
older tools of discrete geometry, such as cuttings.

Let $P$ be a set of $n$ points in the plane, and let $r$ be a
parameter, $1<r\le n$. We say that $f\in\R[x,y]$ is an \notion{$r$-partitioning
polynomial} for $P$ if no connected
component of $\R^2\setminus Z(f)$ contains more than
$n/r$ points of $P$.

In the sequel, we will sometimes
call the connected components of $\R^2\setminus Z(f)$
\emph{cells}. Let us also stress that the cells are open
sets. The points of $P$ lying on $Z(f)$ do not belong to
any cell, and usually they require a special treatment.

\begin{theorem}[Polynomial partitioning theorem]\label{l:ppart}
For every $r>1$, every finite point set
$P\subset\R^2$ admits an $r$-partitioning
polynomial $f$ of degree at most $O(\sqrt{r}\,)$.
\end{theorem}

\heading{Proof.}
We inductively construct collections $\PP_0,\PP_1,\ldots$, each
consisting of disjoint subsets of $P$, such that $|\PP_j|\le 2^j$
for each $j$.  We start with $\PP_0:=\{P\}$. Having constructed
$\PP_j$, with at most $2^j$ sets, we use the polynomial
ham-sandwich theorem to construct a polynomial $f_j$, of degree
$\deg(f_j) \le \sqrt{2\cdot2^j}$, that bisects each of the sets
of $\PP_j$. Then for every subset $Q\in \PP_j$, we let $Q^+$
consist of the points of $Q$ at which $f_j>0$, and let $Q^-$
consist of the points of $Q$ with $f_j<0$, and we put
$\PP_{j+1}:=\bigcup_{Q\in\PP_j}\{Q^+,Q^-\}$.

Each of the sets in $\PP_j$ has size at most $|P|/2^j$. We let $t
=\lceil \log_2 r\rceil$; then each of the sets in $\PP_t$ has
size at most $|P|/r$. We set $f:=f_1f_2\cdots f_t$.

By the construction, no component of $\R^2\setminus Z(f)$ can
contain points of two different sets in $\PP_t$, because any arc
connecting a point in one subset to a point in another subset must
contain a point at which one of the polynomials $f_j$ vanishes, so
the arc must cross $Z(f)$. Thus $f$ is an $r$-partitioning
polynomial for $P$.

It remains to bound the degree:
$$
\deg(f) = \deg(f_1) + \deg(f_2) + \cdots + \deg(f_t) \le \sqrt{2}
\sum_{j=1}^t 2^{j/2} \le \frac{2}{\sqrt{2}-1} 2^{t/2} \le c\sqrt{r} .
$$
where $c = 2\sqrt{2}/(\sqrt{2}-1) < 7$.
\proofend

\heading{A comparison with other partitioning techniques.\footnote{This
part is slightly more advanced and assumes some familiarity
with previous techniques used in incidence problems.} }
The Guth--Katz technique with partitioning polynomials
is useful for problems where we deal with a finite point set $P$
and with a collection $\Gamma$ of lines, algebraic curves,
or algebraic varieties in higher dimensions. It provides
a method of implementing the divide-and-conquer paradigm.

In the planar case discussed above, the plane is
subdivided by $Z(f)$ into some number of (open, connected) cells,
each containing at most $|P|/r$ points of $P$.
If $\Gamma$ consists of lines, then every $\gamma\in\Gamma$
intersects at most $\deg(f)+1=O(\sqrt r)$ cells (by
 Lemma~\ref{l:atmost}).
Similarly, for $\Gamma$ consisting of algebraic curves of
degree bounded by a constant, every $\gamma\in\Gamma$ intersects at most
$O(\sqrt r\,)$ cells by B\'ezout's theorem (Theorem~\ref{t:bezout}).
Thus, if we define,  for every cell $C_i$ of $\R^2\setminus Z(f)$,
a subset $P_i\subseteq P$ as the set of points of $P$ contained in $C_i$,
and we let $\Gamma_i$ consist of the lines or curves of $\Gamma$
intersecting $C_i$, then $|P_i|\le |P|/r$ for all $i$,
 and the \emph{average} size of the $\Gamma_i$ is $O(|\Gamma|/\sqrt r\,)$.

There are two earlier partitioning tools in discrete geometry
with a similar effect. The first, and simpler, kind
of them are \emph{cuttings}
(see \cite{Ch05}). A cutting for a collection
$\Gamma$ of curves in the plane subdivides $\R^2$ into
a collection of connected, simply shaped cells,
in such a way that no cell is crossed
by more than a prescribed fraction of the curves of $\Gamma$.
If we again let $P_i$ denote the set of points of $P$ in the $i$th
cell, and $\Gamma_i$ is the set of the curves intersecting that cell,
then this time \emph{all} $\Gamma_i$ have size $O(|\Gamma|/\sqrt r\,)$,
and the \emph{average} of the sizes $|P_i|$
is\footnote{Here we choose the parameterization so that it agrees
with the one for polynomial partitions; the usual notation
in the literature would use $r$ for our~$\sqrt r$.}
$|P|/r$.  Thus, the behavior of cuttings
is, in a sense, ``dual'' to that of polynomial partitions.
For many applications, this does not really make a difference.

The second of the earlier tools are \emph{simplicial partitions}
\cite{m-ept-92}. Here, as in the case of polynomial
partitions, the plane is subdivided into cells so that
$|P_i|\le |P|/r$ for each $i$
(where, again, $P_i$ is the set of points of $P$ in the $i$th
cell), and no $\gamma\in\Gamma$ intersects more than $O(\sqrt r)$
cells.\footnote{In the original version of simplicial partitions
\cite{m-ept-92}, the cells cover $\R^2$, but they need not
be disjoint. In a newer version due to Chan
\cite{Chan-newpart}, disjointness can also be guaranteed.}

In the plane, as far as we can see, whatever can be done with
polynomial partitions, can also be achieved through cuttings
or through simplicial partitions. The main advantage of polynomial
partitions is simplicity of the proof. On the other hand,
cuttings and simplicial partitions can be constructed and
manipulated with fairly efficient algorithms, at least
in the sense of asymptotic complexity, which is not at all
clear for polynomial partitions. (For example, finding a ham-sandwich
cut in a high-dimensional space is a rather costly operation;
see \cite{Knauer-al} for a computational hardness result
and references.)

Polynomial partitions may be more powerful than the earlier tools
if we pass to a higher-dimensional space $\R^d$, $d>2$.
Asymptotically optimal cuttings
and simplicial partitions are known to exist in $\R^d$,
for every fixed $d$, in the case where $\Gamma$ is a collection
of \emph{hyperplanes}. However, if we want to apply analogous
methods to construct cuttings (or simplicial partitions, whose
construction needs cuttings as a subroutine) for $\Gamma$
consisting of algebraic surfaces of degree bounded by
a constant, say,
then there is a stumbling block. In one of the steps
of the construction, we have
a collection $\Gamma'$ of $m$ surfaces from $\Gamma$. It is known
that these surfaces partition $\R^d$ into $O(m^d)$ cells,
but we need to further subdivide each cell into subcells,
so that each of the resulting subcells can be described by
a constant number of real parameters. There is no known
general solution that achieves $O(m^d)$ subcells in total, which is
the optimal bound one is after for most applications.
For $d=3,4$, the situation is still not bad, since bounds
only slightly worse than $O(m^d)$ have been proved,
but for $d\ge 5$, the best bound is of order roughly $m^{2d-4}$,
and so for  large $d$, the exponent is almost twice larger
of what it probably should be
(see \cite{AM} for a more detailed discussion).
The new approach with polynomial
partitions might hopefully be able to bypass this stumbling block,
at least in non-algorithmic applications.

\section{Proof of the Szemer\'edi--Trotter theorem}\label{s:SzTproof}

We recall that we are given a set $P$ of $m$ distinct points and a set
$L$ of $n$ distinct lines in the plane and we want to bound the
number of incidences $I(P,L)$.

We begin with a simple observation (appearing in most of the previous
proofs).

\begin{lemma} \label{l:kst}
$I(P,L) \le n+m^2$.
\end{lemma}

\heading{Proof.}
We divide the lines of $L$ into two subsets: the lines in
$L'$ are incident to at most one point of $P$,
while the lines in $L''$ pass through at least two points.

Obviously, $I(P,L')\le |L'|\le n$.
In order to bound $I(P,L'')$, we note that a point
$p\in P$ may have at most $m-1$ incidences with  the lines
of $L''$, since there are at most $m-1$ lines passing through $p$
and some other point of~$P$. Thus, $I(P,L'')\le m(m-1)\le m^2$.
\proofend

\medskip
Let us remark that this lemma also follows from
the K\H{o}v\'ari--S\'os--Tur\'an theorem \cite{KST54}
concerning graphs with forbidden complete bipartite subgraphs.
In the above argument, we are really proving
the required instance of K\H{o}v\'ari--S\'os--Tur\'an.
\medskip

\heading{Proof of the
 Szemer\'edi--Trotter theorem. }
For simplicity, we first do the proof for
$m=n$, and then indicate the changes
needed to handle an arbitrary~$m$.

We set $r:=n^{2/3}$, and we let $f$ be an $r$-partitioning
polynomial for $P$. By the polynomial partitioning theorem
(Theorem~\ref{l:ppart}), we may assume $\Deg=\deg(f)=O(\sqrt r\,) =
O(n^{1/3})$.

Let $Z:=Z(f)$, let $C_1,\ldots,C_s$ be the connected components of
$\R^2\setminus Z$,
let $P_i:=P\cap C_i$, and let
$P_0:=P\cap Z$. Since $f$ is an $r$-partitioning polynomial, we have
$|P_i| \le n/r = n^{1/3}$, $i=1,2,\ldots,s$. Furthermore, let
$L_0\subset L$ consist of the lines of $L$ contained in $Z$; we have
$|L_0|\le \Deg$ by Lemma~\ref{l:contfewlines}.

We decompose
$$
I(P,L)= I(P_0,L_0)+I(P_0,L\setminus L_0)
+\sum_{i=1}^s I(P_i,L).
$$
We can immediately bound
$$
I(P_0,L_0)\le |L_0|\cdot |P_0|\le |L_0|n\le \Deg n =O(n^{4/3}),
$$
and
$$
I(P_0,L\setminus L_0)\le |L\setminus L_0|\Deg=O(n^{4/3}),
$$
since each line of $L\setminus L_0$ intersects $Z$, and thus also
$P_0$, in at most $\Deg=\deg(f)$ points.

It remains to bound $\sum_{i=1}^s I(P_i,L)$. Let $L_i\subset L$ be
the set of lines containing at least one point of $P_i$ (the $L_i$ are
typically \emph{not} disjoint). By Lemma \ref{l:kst} we get
$$
\sum_{i=1}^s I(P_i,L_i) \le \sum_{i=1}^s \left(|L_i|+|P_i|^2\right).
$$
We have $\sum_{i=1}^s|L_i|=O((D+1)n)=O(n^{4/3})$, since by
Lemma~\ref{l:atmost}, no line intersects more than $\Deg+1$ of the
sets $P_i$. Finally,
$\sum_{i=1}^s |P_i|^2\le (\max_i |P_i|)\cdot\sum_{i=1}^s|P_i|\le
\frac nr\cdot n=O(n^{4/3})$.
This finishes the proof for the case $m=n$.

We generalize the proof for an arbitrary $m$ as follows. We may assume,
without loss of generality, that $m\le n$;  the complementary case is
handled by interchanging the roles of $P$ and $L$, via a standard
planar duality. We may also assume that $\sqrt{n} \le m$, since
otherwise, the theorem follows from Lemma~\ref{l:kst}.
Then we set $r:=m^{4/3}/n^{2/3}$. Noting that $1\le r\le m$ for the
assumed range of $m$, we then proceed as in the
case $m=n$ above. We get $\Deg=\deg(f) =O(m^{2/3}/n^{1/3})$, and we
check that all the partial bounds in the proof are at most
$O(m^{2/3}n^{2/3})$.
\proofend

\section{Incidences of points with algebraic curves}\label{s:curves}

As was announced in the introduction, we prove the following
(weaker) version of Theorem~\ref{szt:algeb}.

\begin{theorem}\label{t:curveth}
Let $b,k$ and $C$ be constants, let $P$ be a set of $m$ points in the
plane, and let $\Gamma$ be a family of planar curves such that
\begin{enumerate}
\item[\rm(i)] every $\gamma\in\Gamma$ is an algebraic curve of degree
at most~$b$, and
\item[\rm(ii)] for every $k$ distinct points in the plane,
there exist at most $C$ distinct curves in $\Gamma$ passing
through all of them.
\end{enumerate}
Then $I(P,\Gamma)=O\left(m^{k/(2k-1)}n^{(2k-2)/(2k-1)}+m+n\right)$,
with the constant of proportionality depending on $b,k,C$.
\end{theorem}

In the proof, we may assume that the curves in $\Gamma$ are
irreducible.\footnote{We recall that a planar algebraic curve
$\gamma$ is \emph{irreducible} if $\gamma=Z(g)$ for an
irreducible polynomial $g$, i.e., one that cannot be written
as $g=g_1g_2$ with both $g_1,g_2$ nonconstant (and real in our case).
 For $\gamma=Z(g)$
with $g$ arbitrary, we can write $g=g_1g_2\cdots g_k$ as a product
of irreducible factors, and the \emph{irreducible components}
of $g$ are $Z(g_1)$,\ldots, $Z(g_k)$.}
Indeed, if it is not the case, we apply the forthcoming analysis
to the irreducible components of the curves of $\Gamma$,
whose number is at most $bm$.

We begin with an analog of Lemma~\ref{l:kst}.

\begin{lemma} \label{l:kst1} Under the conditions of Theorem~\ref{t:curveth},
we have $I(P,\Gamma) = O(n+m^k)$, and also $I(P,\Gamma)=O(m+n^2)$; the
constants of proportionality depend on $b,k,C$.
\end{lemma}

\heading{Proof.} For the first estimate,
we distinguish between the curves with fewer than $k$
incidences, which altogether generate $O(n)$ incidences, and
curves with at least $k$ incidences, observing that
there are at most $C{m-1\choose k-1}$ such curves through
each point of~$P$.

For the second estimate, we first note that, by the assumed
irreducibility and by B\'ezout's theorem (Theorem~\ref{t:bezout}),
every pair of curves of $\Gamma$ intersect in at most $b^2$ points.
Then we distinguish between points
lying on at most one curve each, which have $O(m)$ incidences
altogether, and the remaining points, each lying on at least
two curves. Now a single $\gamma\in\Gamma$ has at most
$b^2(n-1)$ intersections with the other curves,
and thus it contributes at most $b^2(n-1)$ incidences with these latter points.
So $I(P,\Gamma)=O(m+n^2)$ follows.
\proofend

\heading{Proof of Theorem~\ref{t:curveth}. }
We may assume $m\le n^2$ and $n\le m^{k}$, for otherwise,
the bounds of Lemma~\ref{l:kst1} give $I(P,\Gamma)=O(m+n)$.

We set $r:=m^{2k/(2k-1)}/n^{2/(2k-1)}$, and we observe that
our assumptions on $m,n$ yield $1\le r\le m$.
Let $f$ be an $r$-partitioning polynomial for $P$, of degree
$$
\deg(f) = O(\sqrt{r}\,) = O\bigl(m^{k/(2k-1)}/n^{1/(2k-1)} \bigr) .
$$
The proof now continues in much the same way as the proof of the
Szemer\'edi--Trotter theorem.

We put $Z:= Z(f)$, let $P_0:=P\cap Z$,
and let $\Gamma_0\subset\Gamma$ consist of the curves
fully contained in $Z$. Since every $\gamma\in\Gamma_0$
is irreducible, it must be a zero set of a factor of $f$
(this follows from B\'ezout's theorem), and so
$|\Gamma_0|\le \deg(f)= O(\sqrt r\,)$.
Hence $I(P_0,\Gamma_0)=O(m+|\Gamma_0|^2)= O(m+r)=O(m)$
by the second bound of Lemma~\ref{l:kst1}. (Here
the argument differs from the one for the Szemer\'edi--Trotter
theorem---in the latter, it was sufficient to use the
trivial bound $I(P_0,L_0)\le|P_0|\cdot|L_0|$, which in general is not
sufficient here.)

Next, we consider $\gamma\in\Gamma\setminus\Gamma_0$. Applying B\'ezout's
theorem to $\gamma$ and every irreducible component of $Z$
in turn, we see that $|\gamma\cap Z|\le b\cdot \deg(f)=O(\sqrt r\,)$.
So $I(P_0,\Gamma\setminus \Gamma_0)
=O(n\sqrt r)=O\bigl(m^{k/(2k-1)}
n^{(2k-2)/(2k-1)}\bigr)$.


Letting $C_1,\ldots,C_s$ be the connected
components of $\R^2\setminus Z$,
it remains to bound $\sum_{i=1}^s I(P_i,\Gamma_i)$,
where $P_i=P\cap C_i$ and $\Gamma_i$ is the set of curves
meeting $C_i$. By B\'ezout's theorem once again,
we have $\sum_{i=1}^s|\Gamma_i|=O(n\cdot\deg(f))=
O(n\sqrt r\,)$. Then, by the first bound of Lemma~\ref{l:kst1},
we obtain
\begin{align*}
\sum_{i=1}^s I(P_i,\Gamma_i) &=
O\Bigl(\sum_{i=1}^s \bigl(|\Gamma_i|+|P_i|^k\bigr)\Bigr) \le
O(n\sqrt r\,)+\Bigl(\max{}_i|P_i|\Bigr)^{k-1}
O\Bigl(\sum_{i=1}^s|P_i|\Bigr)\\
&= O(n\sqrt r+(m/r)^{k-1}m)=
O\bigl(m^{k/(2k-1)} n^{(2k-2)/(2k-1)} \bigr).
\end{align*}
\proofend

\section{Spanning trees with low crossing number in the plane}\label{s:pl-sp}

In the forthcoming proof of Theorem~\ref{t:sptree},
instead of constructing a geometric spanning tree directly,
it will be more natural to construct an arcwise connected
set $X$, made of segments and algebraic arcs,
that has a low crossing number and contains the given point set $P$.
Here we say that a set $X\subseteq\R^d$ has \emph{crossing number
at most $k$} if each line, possibly with finitely many exceptions,
intersects $X$ in at most $k$ points. (It is easy to check that
for a geometric spanning tree, this new definition is equivalent
to the earlier one.)

The following  lemma
allows us to convert such an $X$
into a geometric spanning tree. Although we are not aware of
an explicit reference for the statement we need,
most of the ideas of the proof appear in the literature
in some form.

\begin{lemma}\label{l:arcwise2tree}
Let $P$ be a set of $n$ points in the plane, and let
$X$ be an arcwise connected set containing $P$, with crossing
number at most $k$. Then there exists a (geometric) spanning tree
of $P$ whose edges are straight segments and whose crossing number
is at most~$2k$.
\end{lemma}

\heading{Proof.}
In the first stage of the proof we build
a Steiner tree $S$ for $P$, whose edges are arcs contained in $X$.
We order the points of $P$ arbitrarily, into a sequence
$p_1,p_2,\ldots,p_n$. We set $S_1=\{p_1\}$, and, having built a
Steiner tree $S_i\subseteq X$ for $\{p_1,\ldots,p_i\}$, we choose
an arc $\alpha_{i}$ connecting
$p_{i+1}$ to some point $q_i$ of $S_i$, in such a way that
$\alpha_i\cap S_i=\{q_i\}$. Then we set $S_{i+1}:=S_i\cup\alpha_{i}$.
Having reached $i=n$, we set $S:=S_n$; see Fig.~\ref{f:steiner}.
The crossing number of $S$ is at most $k$ since $S\subseteq X$.

\labepsfig{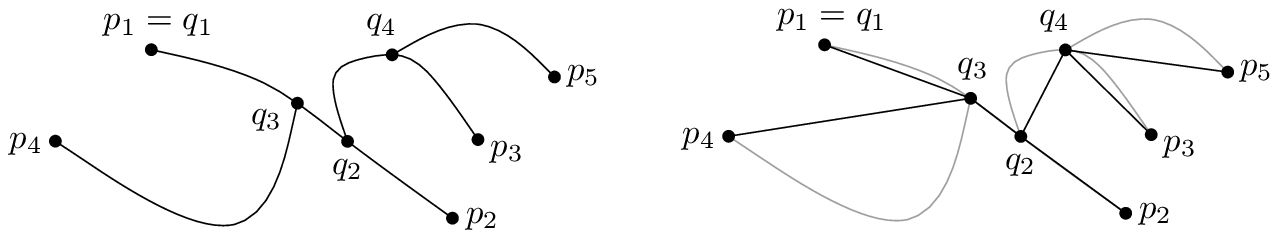}{Illustrating the proof of Lemma~\ref{l:arcwise2tree}:
Left: Building a Steiner tree from arcs. Right: Shortcutting the arcs
into segments.}

In the second stage, we replace arcs by straight segments.
Namely, the points $q_j$ divide $S$ into finitely many
subarcs, and we replace each of them by a straight segment
connecting its endpoints. It is easily seen (and standard)
that the crossing number does not increase.
This yields a Steiner tree for $P$ whose
edges are straight segments.

In the third and last stage, we eliminate the Steiner points
and obtain a spanning tree, at the price of at most doubling
the crossing number. This is done by performing an inorder traversal
of the tree, starting from some arbitrary root vertex,
tracing each edge in both directions, skipping over the
Steiner points, connecting each pair of consecutively visited
points of $P$ by a straight segment, and finally eliminating cycles in
the resulting tour.
\proofend

The main step in the proof of Theorem~\ref{t:sptree} is the following
lemma.

\begin{lemma}\label{l:iteration-sp}
Let $P$ be a set of $n$ points in the plane.
Then there exists a set $X\subseteq \R^2$ that contains $P$,
has at most $n/2$ arcwise connected components, and
with crossing number $O(\sqrt n\,)$.
\end{lemma}

\heading{Proof. } If $n$ is below a suitable constant, we can
interconnect the points of $P$ by an arbitrary geometric spanning tree,
and so we may assume that $n$ is large.

We apply the polynomial partitioning theorem
(Theorem~\ref{l:ppart}), to obtain an $r$-partitioning polynomial
$f$ for $P$, with $r$ as large as possible but so that $Z:=Z(f)$ is
guaranteed to have at most $n/2$ connected components. By
Theorem~\ref{l:ppart}, we have $\deg(f) =O(\sqrt r\,)$, and so, by
Harnack's theorem (Theorem~\ref{t:harnack}), we can afford to take
$r=n/c$ for a suitable constant~$c$.

Then, for every $p\in P$ not lying in $Z$, we pick a straight
segment $\sigma_p$ connecting $p$ to a point of $Z$
(and otherwise avoiding $Z$). We let $X:=Z\cup\bigcup_{p\in P\setminus
Z}\sigma_p$. Clearly, $X$ has at most $n/2$ components,
and it remains to bound its crossing number.

Let $\ell$ be a line that is not contained in $Z$ and that does not
contain any of the segments $\sigma_p$ (these conditions exclude only
finitely many lines). It intersects $Z$ in at most
$\deg(f)=O(\sqrt n\,)$ points, and so it remains to bound the number
of the segments $\sigma_p$ intersected by~$\ell$.

Since $f$ is an $r$-partitioning polynomial for $P$, no component of
$\R^2\setminus Z$ contains more than $c$ points of $P$.  The line
$\ell$ meets at most $1+\deg(f)$ components, and so it intersects at
most $c(1+\deg(f))=O(\sqrt n)$ of the segments $\sigma_p$. The lemma
is proved. \proofend

\heading{Proof of Theorem~\ref{t:sptree}. }
In view of Lemma~\ref{l:arcwise2tree}, it suffices to  construct
an arcwise connected set $X$ containing $P$, with
crossing number $O(\sqrt n)$.
To this end, we apply Lemma~\ref{l:iteration-sp}
recursively.

We construct a sequence $B_0,B_1,B_2,\ldots$ of sets,
such that each $B_i$ contains $P$ and has at most $n/2^{i}$
arcwise connected components. We begin with $B_0:=P$,
and, having constructed $B_i$, we choose a point
in each of its components, which yields a set $R_i$
of at most $n/2^i$ points. Lemma~\ref{l:iteration-sp}
then provides us with a set $X_i\supseteq R_i$
with at most $n/2^{i+1}$ components and crossing
number $O(\sqrt{n/2^i})$. We set $B_{i+1}:=B_i\cup X_i$
and continue with the next iteration, until for some
$i_0$ we reach
an arcwise connected $B_{i_0}$,
which we use as $X$. The crossing numbers of the $X_i$
are bounded by a geometrically decreasing sequence,
and so $X$ has crossing number bounded by its sum, which is
$O(\sqrt n)$, as required.
\proofend

\section{Spanning trees in higher dimensions}\label{s:high-sp}

Here we prove the higher-dimensional
generalization of Theorem~\ref{t:sptree} mentioned in
the introduction.

\begin{theorem}\label{t:highsp}
Every set $P$
of $n$ points in $\R^d$ admits a geometric spanning tree with
crossing number (w.r.t.\ hyperplanes) at most $C_dn^{1-1/d}$,
with $C_d$ a sufficiently large constant depending on~$d$.
\end{theorem}

When one tries to extend the planar proof from Section~\ref{s:pl-sp},
with the appropriate higher-dimensional analogs of the polynomial
partition lemma and Harnack's theorem (discussed below),
a problem arises when
almost all the points of $P$ happen to lie on the zero set $Z(f)$
of the partitioning polynomial. (This situation seems hard to
avoid---for example, $P$ may lie on a low-degree algebraic variety,
in which case the zero set of each of the bisecting polynomials would
simply coincide with this variety.)

In the planar case this did not matter, since we could use $Z(f)$
itself as a part of the connecting set~$X$. However, in higher
dimension, we cannot take all of $Z(f)$ (which is typically
a $(d-1)$-dimensional object), so we would still need to construct
a suitable connecting set with low crossing number within $Z(f)$.

Fortunately, the spanning tree problem behaves well with respect to
small perturbations. Namely, it is easy to see (and well known) that
the crossing number of a geometric spanning tree cannot change by a
(sufficiently small) perturbation of its vertex set, and this will
allow us to avoid the situation with too many points on~$Z(f)$.
Before carrying out this plan, we first summarize and review the
additional tools we need, beyond those already covered.

\subsection{Additional tools}

The polynomial ham sandwich theorem (Theorem~\ref{t:phsc}) and the
polynomial partitioning theorem (Theorem~\ref{l:ppart}) immediately
generalize to $\R^d$, with the $d$-variate $r$-partitioning
polynomial $f$ having degree $O(r^{1/d})$ (this relies on the fact
that the number of monomials of degree $D$ in $d$ variables is
${D+d\choose d}$, so the degree will be the smallest integer
satisfying ${D+d\choose d}-1\ge r$).

We will also need a kind of generalization of Harnack's theorem,
dealing with components of the complement of $Z(f)$, rather than
with the components of~$Z(f)$:

\begin{lemma}\label{l:cpts}
Let $f$ be a real polynomial of degree $D$ in $d$ variables.
Then the number of connected components of
$\R^d\setminus Z(f)$ is at most $6(2D)^d$.
\end{lemma}

This follows, for example, from Warren \cite[Theorem~2]{w-lbanm-68}
(also see \cite{BPR} for an exposition, and
\cite{Akama-al} for a neatly simplified proof).

We also note that if $f$ is as in Lemma~\ref{l:cpts} and $h$ is
a hyperplane in $\R^d$, then $h\setminus Z(f)$ has at most
$6(2D)^{d-1}$ connected components, and consequently,
$h$ intersects at most that many components of $\R^d\setminus Z(f)$.
Indeed, this is clear from Lemma~\ref{l:cpts} if $h$
is the coordinate hyperplane $x_d=0$, and the general
case follows by a linear transformation of coordinates.


\subsection{A general position lemma} \label{sec:per}

We need the following lemma, which is probably known,
but unfortunately we do not have a reference at the moment.

\begin{lemma}\label{l:genpos}
Let $d,\Deg$ be given integers, and let $k:={\Deg+d\choose d}-1$.
Let $P=(p_1,\ldots,p_{k+1})$ be an ordered $(k+1)$-tuple of
points in $\R^d$. Let us call $P$ \emph{exceptional}
if it is contained in the zero set of a
nonzero $d$-variate polynomial of degree at most $\Deg$.
Then there is a nonzero polynomial $\psi=\psi_{d,\Deg}$
with integer coefficients
in the variables $z_{ij}$, $1\le i\le k+1$, $1\le j\le d$,
such that all exceptional $(k+1)$-tuples $(p_1,\ldots,p_{k+1})$
belong to the zero set of $\psi$ (that is,
if we set  $z_{ij}$ to the $j$th coordinate
of $p_i$, for all $i,j$, then $\psi$ evaluates to~$0$).
\end{lemma}

\heading{Proof. }
The value of $k$ in the lemma is the number
of nonconstant monomials of degree at most $\Deg$ in
the $d$ variables $x_1,\ldots,x_d$.
Let $\mu_1,\mu_2,\ldots,\mu_k$ be an enumeration
of these monomials in some fixed order.

It is convenient to phrase the argument
using the Veronese map $\Phi\:\R^d\to \R^k$, which we encountered in
Section~\ref{s:veronese} for the special case $d=2$.
For $x=(x_1,\ldots,x_d)\in\R^d$, we can write
$\Phi(x)=(\mu_i(x) \mid i=1,2,\ldots,k)\in\R^k$.

As in the proof of Theorem~\ref{t:phsc}, the
zero set $Z(f)$ of a polynomial $f$ of degree at most
$\Deg$ can be written as $\Phi^{-1}(h)$, where $h$ is a suitable
hyperplane in $\R^k$. Thus, the condition for a sequence
$P=(p_1,\ldots,p_{k+1})$ of points in $\R^d$ to be exceptional
is equivalent to the $k+1$ points $\Phi(p_1),\ldots,\Phi(p_{k+1})$
lying on a common hyperplane in $\R^k$.

The condition that $k+1$ points $q_1,\ldots,q_{k+1}$
in $\R^k$ lie on a common hyperplane
can be expressed by the vanishing of a suitable determinant
in the coordinates of $q_1,\ldots,q_{k+1}$. Namely,
it is equivalent to $\det(A)=0$, where
$$
A=A(q_1,\ldots,q_{k+1})=\left(\begin{array}{ccccc} 1 & q_{11} & q_{12} & \ldots & q_{1k} \\
1 & q_{21} & q_{22} & \cdots & q_{2k} \\
&\vdots &&\vdots&\\
1& q_{(k+1) 1} & q_{(k+1) 2} & \cdots & q_{(k+1)k} \\
\end{array}
\right).
$$
We define the desired polynomial $\psi$ by
$$\psi=\psi(z_{11},\ldots,z_{(k+1)d}):=
\det(A(\Phi(z_1),\Phi(z_2),\ldots,\Phi(z_{k+1}))),
$$
where $z_i=(z_{i1},\ldots,z_{id})$. Clearly, $\psi$
has integer coefficients,
and by the above, it vanishes on all exceptional sequences;
it remains to verify that it is not identically~0.

Assuming the contrary, it means that the images of any $k+1$ points
under the Veronese map lie on a common hyperplane in $\R^d$. This in
turn implies that all of $\Phi(\R^d)$ is contained in a hyperplane.
But the $\Phi$-preimage of every hyperplane is the zero set of some
nonzero polynomial, and thus it cannot be all of $\R^d$ (extending
the observation in the proof of Lemma~\ref{l:contfewlines} to higher
dimensions).\footnote{Another way to see that $\psi$ is not
identically zero is to consider two terms $t_1$ and $t_2$ in the
expansion of $\det(A(\Phi(z_1),\Phi(z_2),\ldots,\Phi(z_{k+1})))$.
There is a row $i$ of $A(\Phi(z_1),\Phi(z_2),\ldots,\Phi(z_{k+1}))$
from which $t_1$ and $t_2$ contain different elements. It follows
that $t_1$ and $t_2$ must contain a different power of one of the
variables $z_{i1},\ldots,z_{id}$. We get that all terms in the
expansion of $\det(A(\Phi(z_1),\Phi(z_2),\ldots,\Phi(z_{k+1})))$ are
are different monomials, so $\psi$ cannot be identically zero.}
 The resulting
contradiction proves the lemma. \proofend

\subsection{Proof of Theorem~\ref{t:highsp} }

Given a finite
point set $P\subset\R^d$, we first perturb each point slightly,
obtaining a new set $P'$, for which we may assume that the coordinates
of its points are algebraically independent (i.e., they do not satisfy
any nontrivial polynomial equation with integer coefficients).\footnote{The
existence of such $P'$ is well known and follows, e.g., by a standard
measure argument via Sard's theorem (which guarantees that
the zero set of every nonzero multivariate polynomial has zero measure;
see, e.g., \cite{sard,stern}).}

In particular, for every $\Deg=1,2,\ldots$, if we set $k:={\Deg+d\choose d}-1$
as in Lemma~\ref{l:genpos}, then no $(k+1)$-tuple of points
of $P'$ is contained in $Z(f)$, for any nonzero polynomial $f$
of degree at most~$\Deg$.

By the observation mentioned at the beginning of
Section~\ref{s:high-sp}, it suffices to exhibit a geometric spanning
tree with crossing number $O(n^{1-1/d})$ for~$P'$.

Moreover, it suffices to show that there exists a geometric graph
$G$ on the vertex set $P'$ with at most $n/2$ components and with
crossing number $O(n^{1-1/d})$; the existence of the
desired spanning tree then follows
by recursion on the size of $P'$ (as in the proof of the planar
case).

So we set $r:=n/c$ for a sufficiently large constant $c>0$,
and construct an $r$-partitioning polynomial $f$ for $P'$,
of degree $\Deg=O(r^{1/d})$; thus, no component of
$\R^d\setminus Z$ contains more than $c$ points of $P'$,
where $Z=Z(f)$.

By the algebraic independence of $P'$, and by Lemma~\ref{l:genpos},
$Z$ contains fewer than ${\Deg+d\choose d}=O(r)$
points of $P'$, with a constant of proportionality
depending only on $D$. For $c$ sufficiently large, we thus
have $|P'_0|\le \frac n4$, where $P'_0:=P'\cap Z$.
By Lemma~\ref{l:cpts}, we
may also assume that for $c$ sufficiently large, $\R^d\setminus Z$
has at most $\frac n4$ components.

For each component $U$ of $\R^d\setminus Z$, we now interconnect the
points of $U\cap P'$ by an arbitrary geometric spanning tree $T_U$.
The geometric graph $G$ is the union of all the trees $T_U$ and the
points of $P'_0$ (which appear as isolated vertices in $G$).
The number of connected components of $G$ is at most
$\frac n2$ (one for each $T_U$ and one for each point of $P'_0$).
It remains to bound its crossing number.

To this end, let $h$ be a hyperplane avoiding $P'$, and let us
consider an edge $\{p,q\}$ of $G$ crossed by $h$. This edge
belongs to some $T_U$, and so the points $p$ and $q$ lie in
the same component $U$ of $\R^d\setminus Z$. Considering
an arc $\alpha\subset U$ connecting $p$ to $q$, we see that
$h$ has to intersect $\alpha$ and thus $U$ too.
By the remark following Lemma~\ref{l:cpts}, $h$ intersects
at most $O(\Deg^{d-1})= O(n^{1-1/d})$ components of
$\R^d\setminus Z$, and within each such component $U$ it meets at most
$c$ edges of $G$ (that is, of $T_U$). Hence the crossing number of $G$
is $O(n^{1-1/d})$, as claimed.
\proofend

\section{Conclusion}

We regard this paper as an initial stepping stone in the development
of applications of the new algebraic machinery of Guth and Katz.
It is encouraging that this technique can replace more traditional
approaches and yield simpler proofs of central theorems in
combinatorial geometry.

Of course the real challenge is to use the techniques to obtain
improved solutions to other ``hard Erd\H{o}s problems in discrete
geometry'' (borrowing from the title of \cite{s-cnhep-97}), as Guth
and Katz themselves did, first for the joints problem in \cite{GK} and
then for the harder distinct distances problem in \cite{GK2}.
There is a long list of candidate problems, of varying degree of
difficulty. Perhaps the hardest in the list is the planar unit
distances problem of Erd{\H o}s: What is the maximum possible number
of unit distances determined by a set of $n$ points in the plane? This
problem seems to require an algebraic approach, mainly because the
best known upper bound, $O(n^{4/3})$, is known to be tight if the norm
is not Euclidean, as shown by Valtr~\cite{valtr}.

In closing, one should note that the algebraic approach used in this
paper has also some disadvantages. For one, it seems
to require the objects to be algebraic or semialgebraic.
For example, the
Szemer\'edi-Trotter theorem can easily be extended to yield the same
bound on the number of incidences of points and \emph{pseudolines},
using, e.g., the combinatorial proof technique of
Sz\'ekely~\cite{s-cnhep-97}, but such an extension does not seem to
follow from the polynomial partitioning technique. The same situation
occurs in the setup of Theorem~\ref{szt:algeb}, where the general
situation considered there can be handled by traditional combinatorial
tools, but not by the algebraic machinery, which can only establish
weaker variants, like the one in Theorem~\ref{t:curveth}.
Perhaps some abstract version of polynomial
partitions, yet to be discovered, might combine the advantages
of both approaches.

\subsection*{Acknowledgements.}
The authors wish to thank Roel Apfelbaum and Sariel Har-Peled
for useful exchanges of ideas that were helpful in the preparation of
this paper.

\end{document}